\documentclass{article}

\title{Fermat Theorems --- Simple Proofs}
\author{ROBERT J SIBNER}
\usepackage{amsmath}

\begin{document} 

\maketitle
\begin{abstract}
\noindent The Modular Group provides simple proofs of Fermat's representations:\\
\indent$X^2+Y^2$ for primes congruent to 1 (mod 4) and by\\
\indent$X^2+3Y^2$ for primes congruent to 1 (mod 3) 

\end{abstract}
\noindent 1. INTRODUCTION. Among the observations by Fermat [1,4] in the mid seventeenth century about representations of primes as sums of squares were the following:\\[.5pt]

\noindent Theorem A: Primes of the form $p=4N+1$ can be represented as $p=X^2+Y^2$.\\[.5pt]
\noindent Theorem B: Primes of the form $p=3N+1$ can be represented as $p=X^2+3Y^2$.\\[.5pt]

\noindent Since the original proofs by Euler in the following century, many proofs (especially of the first) have been given using a variety of techniques.  However, using the modular group, the arguments that follow provide proofs that are certainly the simplest. In addition, they provide a striking example of the close relation between number theory and complex analysis.  We recall in section 2 some standard facts from algebra and complex analysis and in sections 3 and 4 we present two (essentially equivalent) proofs of each of the theorems; one from the point of view of complex analysis and one in the context of arithmetic groups. Preliminary results were obtained in [2] and [3].
   Proofs of both theorems start, as usual, with preliminary divisibility observations for the relevant prime; for the first that $p=4N+1$ is a factor of $m^2+1$ for some $m$, and for the second that $p=3N+1$ is a factor of $m^2+m+1$ for some $m$.  Following Euler, these facts follow directly from Fermat's Little Theorem by factoring the equation $a^{p-1}-1=0$ ($mod$ $p$) and setting $m=a^N$.\\[4pt]
2. MODULAR GROUP. The modular group (projective special linear group) $PSL(2,Z)$ consists of $2x2$ unimodular matrices with integer entries, where matrices $A$ and $-A$ are identified. An element of this group acts on the upper half complex plane $\mathcal{H}$ as a fractional linear transformation $T(z)=(az+b)(cz+d)^{-1}$ (with determinant $ad-bc=1$); the group $\Gamma$ of transformations also being referred to as the modular group. It is well known (and easy to verify) that the set \{${z\in \mathcal{H}: |z|>1, |Re(z)|<1/2}$ \}, together with that part of its boundary with $Re(z)\geq 0$ is a \textit{fundamental set} $\mathcal{F}$; every point in $\mathcal{H}$ is $\Gamma-$equivalent to a point in $\mathcal {F}$ and no two points of $\mathcal{F}$ are $\Gamma-$equivalent.

\noindent Among the elements of the modular group are the elliptic elements of order two ($E_{(2)}^{2} =-$identity).  They have zero trace and, with $r=(1+m^2)s^{-1}$, the form
\begin{equation}
E_{(2)}=\left( \
\begin{array}{cc}
m&-r\\
s&-m
\end{array} 
\right)
\end{equation} 
Elements of the two conjugacy classes of elliptic elements of order three \\    ($E_{(3)}^{3}=-$identity) have trace one and the forms (with $r=(1+m+m^2)s^{-1}$)
\begin{equation}
E_{(3)}=\left( \
\begin{array}{cc}
1+m&-r\\
s&-m
\end{array} 
\right)
\end{equation} and \begin{equation}
E_{(3)}^{-1}=\left( \
\begin{array}{cc}
-m&-r\\
s&1+m
\end{array} 
\right)
\end{equation}
Note that, if $s=1$ or a prime of the form $p=4N+1$ in the first case ($3N+1$ in the second or third), then $r$ is an integer and the matrix is in $PSL(2,Z)$.  As elements of $\Gamma, E_{(2)}$ has the "elliptic" fixed point $(m+i)/p$ while both $E_{(3)}$ and $E_{(3)}^{-1}$ have the fixed point  $[2m+1+i\sqrt{3}]/2p.$\\[4pt]
\noindent 3.1 PROOF OF THEOREM A: The fixed point $(m+i)/p$ of the transformation corresponding to $E_{(2)}$(with $s$ the prime $p=4N+1)$ must be $\Gamma-$ equivalent, by a transformation $T$ in $\Gamma$, to the fixed point $z=i$ of the transformation $z\to-1/z$, the unique order two elliptic fixed point in the fundamental set $\mathcal{F}$ of $\Gamma$.  Comparing the imaginary parts of $T(i)=(a+ib)/(c+id)=(bd-ac+i)/(c^2+d^2)$ with that of $(m+i)/p$ gives immediately $p=c^2+d^2$ \\

\noindent 3.2 PROOF OF THEOREM B. With $s$ the prime $p=3N+1$, both $E_{(3)}^{-1}$ and $E_{(3)}$ have the fixed point $[2m+1+i\sqrt{3}]/2p$ which must be equivalent to the unique fixed point $[1+i\sqrt{3}]/2$ of order 3 in $\mathcal{F}$.  Comparing imaginary parts of $[2m+1+i\sqrt{3}]/2p$ with $T([1+i\sqrt{3}]/2)$ one finds, at first, that $p=d^2+dc+c^2$. If either $c$ or $d$ are even, we can write $p=(d+c/2)^2+3(c/2)^2$ or $p=(c+d/2)^2+3(d/2)^2$ while if both $c$ and $d$ are odd, $p= ([d-c]/2)^2 + 3([d+c]/2)^2$.  In any event,  as claimed in Fermat Theorem B, $p$ is representable by the form $X^2+3Y^2$.\\
\pagebreak

\noindent 4: VARIATION OF PROOFS. The above simple arguments can be translated into equally simple ones about the arithmetic group $PSL(2,Z)$.  In this group there is only one conjugacy class of order two elliptic elements so the matrix $E_{(2)}$ with $p=4N+1$  (and corresponding $m$) is conjugate to the matrix with $s=1$ and $m=0$.  Performing the conjugation by a matrix 
$\left(
\begin{smallmatrix}
a & b\\
c & d
\end{smallmatrix}
\right)$
establishes that $p=c^2+d^2$ proving Fermat's Theorem A.

\noindent For a given $p$ the matrices $E_{(3)}$ and $E_{(3)}^{-1}$ of section 2 are conjugate respectively to the corresponding matrices with $s=1$ and $m=0$.  Performing the conjugations and comparing the result, one obtains in either case that $c^2+dc+d^2$. As before, one sees that this gives the representation $p=X^2+3Y^2$ of Theorem B.\\[6pt]

				REFERENCES\\[4pt]
\noindent 1. FERMAT, Oeuvres (standard Tannery Henry edition) vol. II pp. 113,403.\\
2. K. HASHIMOTO and R. J. SIBNER, Involutive Modular Transformations on the Siegel Upper Half Space and an Application to Representations of Quadratic Forms, J. Number Theory 23, ((1986) 102-110.\\
3 R. J. SIBNER, The extended modular group and sums of squares, Complex Variables Theory Appl. 3 (1984), 481-491.\\
4. A. WEIL, Number Theory, Birkhauser, Boston (1984).\\[7pt]
\noindent CITY UNIVERSITY OF NEW YORK, Graduate Center and Brooklyn College\\*
\nopagebreak[4]Email address: rsibner@gmail.com, rsibner@gc.cuny.edu

\end{document}